\documentclass{amsart}
\usepackage{amssymb}
\usepackage{amsmath}
\usepackage{amsfonts}

\newtheorem{theorem}{Theorem}
\theoremstyle{plain}

\newtheorem{example}{Example}

\newtheorem{lemma}{Lemma}

\newtheorem{remark}{Remark}

\numberwithin{equation}{section}
\input{tcilatex}

\begin{document}
\title{A CLASS OF 3-DIMENSIONAL CONTACT METRIC MANIFOLDS}
\author{I. K\"{u}peli Erken}
\address{Art and Science Faculty,Department of Mathematics, Uludag
University, 16059 Bursa, TURKEY}
\email{iremkupeli@uludag.edu.tr}
\author{C. Murathan}
\address{Art and Science Faculty,Department of Mathematics, Uludag
University, 16059 Bursa, TURKEY}
\email{cengiz@uludag.edu.tr}
\date{18.05.2012}

\begin{abstract}
We classify the contact metric $3$-manifolds that satisfy $\left\Vert
grad\lambda \right\Vert =1$ and $\bigtriangledown _{\xi }\tau =2a\tau \phi .$
\end{abstract}

\maketitle

\section{\textbf{Introduction}}

Let $M$ be a Riemannian manifold. The tangent sphere bundle $T_{1}M$ admits
a contact metric structure $(\phi ,\xi ,\eta ,g)$ and so $T_{1}M$\ together
with this structure is a contact metric manifold \cite{BL3}. If $M$ is of
constant sectional curvature, then the curvature tensor $R$ of $T_{1}M$ $%
(\phi ,\xi ,\eta ,g)$ satisfies the condition%
\begin{equation}
R(X,Y)\xi =\kappa (\eta (Y)X-\eta (X)Y)+\mu (\eta (Y)hX-\eta (X)hY)
\label{*}
\end{equation}%
for any $X,Y\in \chi (T_{1}M)$, where $2h$ is the Lie derivative of $\phi $
with respect to $\xi $ and $\kappa $,$\mu $ are constants. The contact
metric manifolds satisfying this condition are called $(\kappa ,\mu )$%
-contact metric manifolds. This kind of manifold was introduced and deeply
studied by Blair, Koufogiorgos and Papantoniou in \cite{BL1}. In \cite{FP1}, 
\cite{FP2} and \cite{S} are studied contact metric manifolds satisfying (\ref%
{*}) but with $\kappa ,\mu $ \ smooth functions not necessarily constant. If
\ $\kappa ,\mu $ are nonconstant smooth functions on $M$, the manifold $M$
is called a generalized ($\kappa ,\mu )$-contact metric manifold (see \cite%
{KOU1}). Koufogiorgos \ gave the examples satisfying (\ref{*}) with $\kappa
,\mu $ \ non constant smooth functions for dimension $3$ and also proved in 
\cite{KOU1} that if $\dim M>3$ then $\kappa $,$\mu $ were necessarely
constant. In \cite{KOU2} Koufogiorgos, Markellos and Papantoniou proved the
existence of a new class of contact metric manifolds: the so called $(\kappa
,\mu ,v)$-contact metric manifolds. Such a manifold $M$ is defined through
the condition%
\begin{equation}
R(X,Y)\xi =\kappa (\eta (Y)X-\eta (X)Y)+\mu (\eta (Y)hX-\eta (X)hY)+\nu
(\eta (Y)\phi hX-\eta (X)\phi hY)  \label{kmv}
\end{equation}%
for every $X,Y\in \chi (M)$ and $\kappa ,\mu ,\nu $ are smooth functions on $%
M.$ Furthermore, it is shown in \cite{KOU2}\ \ that if dim $M>3,$ then $%
\kappa ,\mu $ are constants and $\nu $ is the zero function. In \cite{KOU3}
locally classified the $3$-dimensional generalized $(\kappa ,\mu )$-contact
metric manifolds satisfied the condition $\left\Vert grad\kappa \right\Vert
= $constant$(\neq 0)$ . In \cite{KOU4} locally classified the $3$%
-dimensional $(\kappa ,\mu ,v)$-contact metric manifolds which satisfy the
condition $\nabla _{\xi }\tau =2a\tau \phi ,$ where $a$ is a smooth function.

In this paper we want to generalize \ the corollaries of \cite{KOU3} and 
\cite{KOU4} . The paper is organized in the following way. The Section $2$
contains the presentation of some basic notions about contact manifolds and $%
(\kappa ,\mu )$-contact metric manifolds, $(\kappa ,\mu ,v)$-contact metric
manifolds. In section $3$ we give some properties of contact metric $3$%
-dimensional manifold. In section $4$ we generalize the corollaries of \cite%
{KOU3} and \cite{KOU4}. Finally we gave two examples which satisfiy the
conditions of this manifold.

\section{Preliminaries}

A differentiable manifold $M$ of dimension $2n+1$ is said to be a contact
manifold if it carries a global 1-form $\eta $ such that $\eta \wedge (d\eta
)^{n}\neq 0$. It is well known that then there exists a unique vector field $%
\xi $ (called the Reeb vector field) such that $\eta (\xi )=1$ and $d\eta
(\xi ,\cdot )=0$. It is well known that there also exists a Riemannian
metric $g$ and a $(1,1)$-tensor field $\phi $ such that 
\begin{gather}
\phi (\xi )=0,\text{ \ \ }\phi ^{2}=-I+\eta \otimes \xi ,\ \ \eta \circ \phi
=0,  \label{B1} \\
~g(\phi X,\phi Y)=g(X,Y)-\eta (X)\eta (Y)  \label{B0}
\end{gather}%
for any vector field $X$ and $Y$ on $M$. The structure $(\phi ,\xi ,\eta ,g)$
can be chosen so that $d\eta (X,Y)=g(X,\phi Y).$ The manifold $M$ together
with the structure tensors $(\phi ,\xi ,\eta ,g)$ is called a contact metric
manifold structure and is denoted by $M(\phi ,\xi ,\eta ,g)$. Define an
operator $h$ by $h=\frac{1}{2}{\mathcal{L}}_{\xi }\phi $, where $\mathcal{L}$
denotes Lie differentiation. The tensor field $h$ vanishes identically if
and only if the vector field $\xi $ is Killing and in this case the contact
metric manifold is said to be \textit{K-contact}. It is well known that $h$
and $\phi h$ are symmetric operators, $h$ anti-commutes with $\phi $ 
\begin{equation}
\phi h+h\phi =0,\text{ }h\xi =0,\text{ }\eta \circ h=0,\text{ }trh=\text{ }%
tr\phi h=0,  \label{B1a}
\end{equation}%
where $trh$ denotes the trace of $h$. Since $h$ anti-commutes with $\phi $,
if $X$ is an eigenvector of $h$ corresponding to the eigenvalue $\lambda $
then $\phi X$ is also an eigenvector of $h$ corresponding to the eigenvalue $%
-\lambda $ . Moreover, for any contact manifold $M$, the following is
satisfied%
\begin{equation}
\nabla _{X}\xi =-\phi X-\phi hX  \label{B2}
\end{equation}%
where $\nabla $ is the Riemannian connection of $g$.

\textbf{\ }On a contact metric manifold $M^{2n+1}$ we have the formulas

\begin{eqnarray}
(\nabla _{\xi }h) &=&\phi (I-h^{2}-l),  \label{B3} \\
l-\phi l\phi &=&-2(h^{2}+\phi ^{2}),  \label{B4} \\
Trl &=&g(Q\xi ,\xi )=2n-Trh^{2},  \label{B5} \\
\tau &=&2g(\phi \cdot ,h\cdot ),  \label{B6} \\
\nabla _{\xi }\tau &=&2g(\phi \cdot ,\nabla _{\xi }h\cdot ),  \label{B7} \\
\left\Vert \tau \right\Vert ^{2} &=&4\text{tr}h^{2}  \label{B8}
\end{eqnarray}%
where $l=R(X,\xi )\xi $,$~Q$ is Ricci operator of $M.$

A contact metric manifold satisfying $R(X,Y)\xi =0$ is locally isometric to $%
E^{n+1}\times S^{n}(4)$ for $n>1$ and flat for $n=1$ \cite{BL2}.

If a contact metric manifold $M$ is normal (i.e., $N_{\phi }+2d\eta \otimes
\xi =0,$ where $N_{\phi }$ denotes the \textit{Nijenhuis tensor} formed with 
$\phi $), then $M$ is called a \textit{Sasakian manifold}. Equivalently, a
contact metric manifold is Sasakian if and only if $\ (\nabla _{X}\phi
)Y=g(X,Y)\xi -\eta (Y)X$ $\ $or $\ R(X,Y)\xi =\eta (Y)X-\eta (X)Y$ \cite{BL3}%
.

As a generalization of both $R(X,Y)\xi =0$ and the Sasakian manifold consider

\begin{equation}
R(X,Y)\xi =\kappa (\eta (Y)X-\eta (X)Y)+\mu (\eta (Y)hX-\eta (X)hY)
\label{nullity}
\end{equation}%
for smooth functions $\kappa $ and $\mu $ on $M$. If \ $\kappa ,\mu $ is
constant $M$ is called $(\kappa ,\mu )$ contact metric manifold. Otherwise $%
M $ is generalized ($\kappa ,\mu )$ contact metric manifold. Since $15$
years, this kind of manifolds especially was introduced and studied by
Blair, Koufogiorgos, Papantoniou and Markellos in \cite{BL1},\cite{KOU2}.

Let $M(\phi ,\xi ,\eta ,g)$ be a contact metric manifold.$~$A D-homothetic
transformation\textbf{\ }\cite{Tanno} is the transformation

\begin{equation}
\bar{\eta}=\alpha \eta ,\text{ \ \ }\bar{\xi}=\frac{1}{\alpha }\xi ,\text{ \
\ }\bar{\phi}=\phi ,\text{ \ \ }\bar{g}=\alpha g+\alpha (\alpha -1)\eta
\otimes \eta  \label{D}
\end{equation}%
at the structure tensors, where $\alpha $ is a positive constant. It is well
known \cite{Tanno} that $M(\bar{\phi},\bar{\xi},\bar{\eta},\bar{g})$ is also
a contact metric manifold. When two contact structures ($\phi ,\xi ,\eta ,g)$
and $(\bar{\phi},\bar{\xi},\bar{\eta},\bar{g})$ are related by (\ref{D}) ,
we will say that they are D-homothetic\textbf{.}

We can easily show that \ $\bar{h}=\frac{1}{\alpha }h$ so $\bar{\lambda}=%
\frac{1}{\alpha }\lambda .$

Using the relations above we finally obtain that

\begin{equation*}
\bar{R}(X,Y)\bar{\xi}=\frac{\kappa +\alpha ^{2}-1}{\alpha ^{2}}(\bar{\eta}%
(Y)X-\bar{\eta}(X)Y)+\frac{\mu +2(\alpha -1)}{\alpha }(\bar{\eta}(Y)\bar{h}X-%
\bar{\eta}(X)\bar{h}Y).
\end{equation*}%
for all vector fields $X$ and $Y$ on $M.$ Thus $M(\bar{\phi},\bar{\xi},\bar{%
\eta},\bar{g})\ $is a\textit{\ }$(\bar{\kappa},$\ $\bar{\mu})-$contact
metric manifold with

\begin{equation*}
\bar{\kappa}=\frac{\kappa +\alpha ^{2}-1}{\alpha ^{2}},\text{ \ \ }\bar{\mu}=%
\frac{\mu +2(\alpha -1)}{\alpha }.
\end{equation*}%
It is well known, for example, \cite{BL2} , that every $3$-dimensional
contact metric manifold satisfies the integrability condition,

\begin{equation*}
(\nabla _{X}\phi )Y=g(X+hX,Y)\xi -\eta (Y)(X+hX).
\end{equation*}%
Now we will introduce a generalized $(\kappa ,\mu )$-contact metric manifold.

\begin{example}
\cite{KOU1, KOU4} We consider the 3-dimensional manifold $M=\left\{
(x_{1,}x_{2},x_{3})\in R^{3}\mid x_{3}\neq 0\right\} ,$ where $%
(x_{1,}x_{2},x_{3})$ are the standart coordinates in $R^{3}.$ The vector
fields%
\begin{equation*}
e_{1}=\frac{\partial }{\partial x_{1}},\ e_{2}=-2x_{2}x_{3}\frac{\partial }{%
\partial x_{1}}+\frac{2x_{1}}{x_{3}^{3}}\frac{\partial }{\partial x_{2}}-%
\frac{1}{x_{3}^{2}}\frac{\partial }{\partial x_{3}},\ e_{3}=\frac{1}{x_{3}}%
\frac{\partial }{\partial x_{2}}
\end{equation*}

are linearly independent at each point of $M$. Let $g$ be the Riemannian
metric defined by $g(e_{i},e_{j})=\delta _{ij},$ $i,j=1,2,3$ \ $\ $and $\eta 
$ the dual $1$-form to the vector field $e_{1}.$We define the tensor $\phi $
of type $(1,1)$ by $\phi e_{1}=0,\phi e_{2}=e_{3},\phi e_{3}=-\ e_{2}.$%
Following \cite{KOU1} , we have that $M(\eta ,e_{1},\phi ,g)$ is a
generalized $(\kappa ,\mu )$-contact metric manifold with $\kappa =\frac{%
x_{3}^{4}-1}{x_{3}^{4}},$ $\mu =2\left( 1-\frac{1}{x_{3}^{2}}\right) $.By a
straightforward calculation, one can deduce that $M,$ satisfies%
\begin{equation*}
\nabla _{\xi }\tau =2\left( 1-\frac{1}{x_{3}^{2}}\right) \tau \phi .
\end{equation*}
\end{example}

In \cite{KOU2} the authors proved the existence of a new class of contact
metric manifolds which is called $(\kappa ,\mu ,\upsilon )$-contact metric
manifold. This means that curvature tensor $R$ satisfies the condition%
\begin{eqnarray*}
~R(X,Y)\xi &=&\kappa (\eta (Y)X-\eta (X)Y)+\mu (\eta (Y)hX-\eta (X)hY) \\
&&+\upsilon (\eta (Y)\phi hX-\eta (X)\phi hY)
\end{eqnarray*}%
for any vector fields $X$ , $Y$ and $\kappa ,\mu ,\upsilon $ are smooth
functions.

Furthermore, it is shown in \cite{KOU2} that if $dimM>3$, then $\kappa ,\mu $
are constants and $\upsilon $ is the zero function.

\bigskip

\begin{example}
\cite{KOU2} Let $M$ be 3-dimensional contact metric manifold such that%
\begin{equation*}
M=\left\{ (x_{,}y,z)\in R^{3}\mid x>0,\text{ }y>0,\text{ }z>0\right\} ,
\end{equation*}

where ($x,y,z)$ are the cartesian coordinates in $R^{3}.$ We define three
vector fields on $M$ as
\end{example}

\begin{equation*}
e_{1}=\frac{\partial }{\partial x},\ e_{2}=\frac{\partial }{\partial y},\ \
e_{3}=-\frac{4}{z}e^{G}\ G_{y}\frac{\partial }{\partial x}+\beta \ \frac{%
\partial }{\partial y}+e^{G/2}\ \ \frac{\partial }{\partial z}
\end{equation*}%
\textit{where }$G=G(y,z)<0$\textit{\ for all }$(y,z)$\textit{\ is a solution
of the partial differential equation}

\begin{equation*}
2G_{yy}+G_{y}^{2}=-ze^{-G},
\end{equation*}%
\textit{and the function }$\beta =\beta (x,y,z)$\textit{\ solves the system
of partial differential equations}

\begin{eqnarray*}
\beta _{x} &=&\frac{4}{zx^{2}}e^{G}, \\
\beta _{y} &=&\frac{1}{2z}e^{G/2}-\frac{G_{z}e^{G/2}}{2}-\frac{4e^{G}G_{y}}{%
xz}.
\end{eqnarray*}%
\textit{Setting }$\kappa =1-(4e^{2G)}/(z^{2}x^{4}),$\textit{\ \ \ }$\mu
=2(1+(2e^{G})/(zx^{2})),$\textit{\ and }$\nu =-2/x.$

\textit{By direct calculation, these relations yield}%
\begin{eqnarray*}
R(Z,W)\xi &=&\kappa (\eta (W)Z-\eta (Z)W)+\mu (\eta (W)hZ-\eta (Z)hW)+ \\
&&+\nu (\eta (W)\phi hZ-\eta (Z)\phi hW)
\end{eqnarray*}%
\textit{for all vector fields }$Z,W,$\textit{\ on }$M,$\textit{\ where }$%
\kappa ,\mu ,\nu $\textit{\ are nonconstant smooth functions. Hence, it has
been shown that }$M$\textit{\ is a (generalized) }$(\kappa ,\mu ,\nu )$%
\textit{-contact metric manifold.}

\section{Three dimensional contact metric manifolds}

In this section, we will give some properties of contact metric $3-$%
dimensional manifold.

Let $M(\phi ,\xi ,\eta ,g)$ be a contact metric $3$-manifold. Let%
\begin{eqnarray*}
U &=&\left\{ p\in M\mid h(p)\neq 0\right\} \subset M \\
U_{0} &=&\left\{ p\in M\mid h(p)=0\right\} \subset M
\end{eqnarray*}

That $h$ is a smooth function on $M$ implies $U\cup U_{0}$ is an open and
dense subset of $M$, so any property satisfied in $U_{0}\cup U$ is also
satisfied in $M.$

For any point $p\in U\cup U_{0}$, there exists a local orthonormal basis $%
\left\{ e,\phi e,\xi \right\} $ of smooth eigenvectors of h in a
neighbourhood of $p~$(this we call a $\phi $-basis).

On $U$, we put $he=\lambda e,$ $h\phi e=-\lambda \phi e,$where $\lambda $ is
a nonvanishing smooth function assumed to be positive.

\begin{lemma}
\cite{FP1} On the open set $U$ we have%
\begin{eqnarray}
\nabla _{\xi }e &=&a\phi e,\text{ }\nabla _{e}e=b\phi e,\text{ }\nabla
_{\phi e}e=-c\phi e+(\lambda -1)\xi ,  \label{3.1} \\
\nabla _{\xi }\phi e &=&-ae,\nabla _{e}\phi e=-be+(1+\lambda )\xi ,\nabla
_{\phi e}\phi e=ce,  \label{3.2} \\
\nabla _{\xi }\xi &=&0,\text{\ }\nabla _{e}\xi =-(1+\lambda )\phi e,\text{ }%
\nabla _{\phi e}\xi =(1-\lambda )e,  \label{3.3} \\
\nabla _{\xi }h &=&-2ah\phi +(\xi \cdot \lambda )s  \label{3.4}
\end{eqnarray}
\end{lemma}

\textit{where }$a$\textit{\ is a smooth function,}

\begin{eqnarray}
b &=&\frac{1}{2\lambda }((\phi e\cdot \lambda )+A)\text{ \ with \ }A=\eta
(Qe)=S(\xi ,e),  \label{3.5} \\
c &=&\frac{1}{2\lambda }((e\cdot \lambda )+B)\text{ \ with \ }B=\eta (Q\phi
e)=S(\xi ,\phi e),  \label{3.6}
\end{eqnarray}

\textit{and }$s$\textit{\ is the type }$(1,1)$\textit{\ tensor field defined
by }$s\xi =0$\textit{, }$se=e$\textit{\ and }$s\phi e=-\phi e.$

\textit{From Lemma and the formula }$\left[ X,Y\right] =\nabla _{X}Y-\nabla
_{Y}X$\textit{, we can prove that}%
\begin{eqnarray}
\left[ e,\phi e\right] &=&\nabla _{e}\phi e-\text{\ }\nabla _{\phi
e}e=-be+c\phi e+2\xi ,  \label{20} \\
\left[ e,\xi \right] &=&\nabla _{e}\xi -\nabla _{\xi }e=-(a+\lambda +1)\phi
e,  \label{21} \\
\left[ \phi e,\xi \right] &=&\nabla _{\phi e}\xi -\nabla _{\xi }\phi
e=(a-\lambda +1)e.  \label{22}
\end{eqnarray}

\section{Main Results}

First of all, we will introduce the theorem which \ gives motivation to us.

In , \cite{KOU3} locally classified the $3$-dimensional generalized $(\kappa
,\mu )$-contact metric manifolds satisfied the condition $\parallel
grad\kappa \parallel =$const.($\neq 0).$ This classification is given
following.

\begin{theorem}
\cite{KOU3} Let $M(\phi ,\xi ,\eta ,g)$ be a $3$-dimensional generalized $%
(\kappa ,\mu )$-contact metric manifold with $\parallel grad\kappa \parallel
=1.$ Then at any point $p\in M$ there exist a chart $(U,(x,y,z)),z<1$, such
that $\kappa =z$ and $\mu =2(1-\sqrt{1-z)\text{ }}$or $\mu =2(1+\sqrt{1-z})$%
. In the first case ($\mu =2(1-\sqrt{1-z\text{ }}))$, the following are
valid,
\end{theorem}

\begin{equation*}
\xi =\frac{\partial }{\partial x},~~\phi X=\frac{\partial }{\partial y}~~%
\text{and }X=a\frac{\partial }{\partial x}+b\frac{\partial }{\partial y}+%
\frac{\partial }{\partial z}.
\end{equation*}

\textit{In the second case (}$\mu =2(1+\sqrt{1-z})),$\textit{\ the following
are valid,}

\begin{equation*}
\xi =\frac{\partial }{\partial x},\mathit{~~}X=\frac{\partial }{\partial y}%
\mathit{~~}\text{and }\phi X=a^{\prime }\frac{\partial }{\partial x}%
+b^{\prime }\frac{\partial }{\partial y}+\frac{\partial }{\partial z},
\end{equation*}

\textit{where }$a(x,y,z)=-2y+f(z),$ $a^{\prime
}(x,y,z)=2y+h(z),~b(x,y,z)=b^{\prime }(x,y,z)=2x\sqrt{1-z}+\frac{y}{4(1-z)}%
+r(z)$\textit{\ and }$f,r,h$\textit{\ are smooth functions of }$z$.

In \cite{KOU4} the authors \ proved following theorem.

\begin{theorem}
\cite{KOU4} Let $M(\phi ,\xi ,\eta ,g)$ be 3-dimensional $(\kappa ,\mu
,\upsilon )$-contact metric manifold for which $\nabla _{\xi }\tau =2a\tau
\phi $ where $a$ is a smooth function on $M.$ Then $M$ is either a Sasakian
manifold $(\kappa =1)$ or a generalized $(\kappa ,\mu )$-contact metric
manifold with $\kappa <1.$
\end{theorem}

Moreover, the authors locally classified three dimensional $(\kappa ,\mu
,\upsilon )$-contact metric manifolds which satisfy the condition $\nabla
_{\xi }\tau =2a\tau \phi $ in \cite{KOU4}. They gave the following theorem.

\begin{theorem}
\cite{KOU4} Let $M(\phi ,\xi ,\eta ,g)$ be a non-Sasakian $3$-dimensional ($%
\kappa ,\mu ,\upsilon )$-contact metric manifold. Suppose that $\nabla _{\xi
}\tau =2a\tau \phi $ where $a$ is a smooth function which is constant along
the geodesic foliation generated by $\xi .$ Then, $M$ is a generalized $%
(\kappa ,\mu )$-contact metric manifold with $\xi (\mu )=0$. We denote by $%
\left\{ \xi ,e,\phi e\right\} $ a local orthonormal frame of eigenvectors of 
$h$ such that $he=\lambda e$, $\lambda =\sqrt{1-\kappa }>0$. Furthermore,
for every $p\in M$ there exists a chart $(U,(x,y,z)),z<1,$ such that the
function $\kappa $ depends only of the variable $z$ and $\mu =2(1-\sqrt{%
1-\kappa })$ or $\mu =2(1+\sqrt{1-\kappa }).$

In the first case, the following are valid,%
\begin{equation*}
\xi =\frac{\partial }{\partial x},~~\phi e=\frac{\partial }{\partial y}~~%
\text{and }e=\alpha \frac{\partial }{\partial x}+b\frac{\partial }{\partial y%
}+\frac{\partial }{\partial z}.
\end{equation*}

In the second case, the following are valid,%
\begin{equation*}
\xi =\frac{\partial }{\partial x},e=\frac{\partial }{\partial y}~\text{and}%
~~\phi e=\alpha _{1}\frac{\partial }{\partial x}+b_{1}\frac{\partial }{%
\partial y}+\frac{\partial }{\partial z},
\end{equation*}

where 
\begin{eqnarray*}
\alpha (x,y,z) &=&-2y+f(z),\alpha _{1}(x,y,z)=2y+f(z) \\
b(x,y,z) &=&b_{1}(x,y,z)=2\lambda (z)x-\frac{\lambda ^{\prime }(z)}{2\lambda
(z)}y+\psi (z) \\
\lambda (z) &=&\sqrt{1-\kappa (z)}>0
\end{eqnarray*}

and $f,\psi $ are smooth functions of $z.$

\begin{remark}
(a) As a result, contact metric manifold with $\left\Vert grad\text{ }%
\lambda \right\Vert _{g}=d\neq 0$ (cons.) is $D_{\alpha }$- deformed in
another contact metric manifold. with $\left\Vert grad\text{ }\bar{\lambda}%
\right\Vert _{\bar{g}}=d\alpha ^{-\frac{3}{2}}$ and choosing $\alpha =d^{%
\frac{2}{3}},$ it is enough to study those contact metric manifold. with $%
\left\Vert grad\text{ }\lambda \right\Vert =1$.

(b)If $d=0,~$then $\lambda $ is constant. As a result, if $\lambda =0,$ then 
$M$ is a Sasakian manifold.
\end{remark}
\end{theorem}

Now we will give our main theorem.

\begin{theorem}
Let $M(\phi ,\xi ,\eta ,g)$ be a $3$-dimensional contact metric manifold
with $\left\Vert grad\text{ }\lambda \right\Vert =1$ and $\nabla _{\xi }\tau
=2a\tau \phi .$ Then at any point $p\in M$ there exist a chart $(U,(x,y,z))$
such that $\lambda =g(z)\neq 0$ and $A=0,B=F(y,z)$ or $A=F(y,z),B=0.$ In the
first case ($A=0,B=F(y,z)$), the following are valid,

\begin{equation*}
\xi =\frac{\partial }{\partial x},~~\phi e=\frac{\partial }{\partial y}~~%
\text{and }e=k_{1}\frac{\partial }{\partial x}+k_{2}\frac{\partial }{%
\partial y}+k_{3}\frac{\partial }{\partial z},\text{ \ }k_{3}\neq 0.
\end{equation*}%
In the second case ($A=F(y,z),B=0$), the following are valid,

\begin{equation*}
\xi =\frac{\partial }{\partial x},\text{ }e=\frac{\partial }{\partial y}~%
\text{and}~~\phi e=k_{1}^{\prime }\frac{\partial }{\partial x}+k_{2}^{\prime
}\frac{\partial }{\partial y}+k_{3}^{\prime }\frac{\partial }{\partial z},%
\text{ \ \ }k_{3}^{\prime }\neq 0~,
\end{equation*}%
where, 
\begin{equation*}
k_{1}(x,y,z)=-2y+r(z),\text{ \ }k_{1}^{\prime }(x,y,z)=2y+r^{\prime }(z),
\end{equation*}%
\begin{equation*}
k_{2}(x,y,z)=k_{2}^{\prime }(x,y,z)=2xg(z)-\frac{(H(y,z)+y)}{2g(z)}+\beta
(z),
\end{equation*}%
\begin{equation*}
k_{3}(x,y,z)=k_{3}^{\prime }(x,y,z)=t(z)+\delta ,\text{ \ \ }\frac{\partial
H(y,z)}{\partial y}=F(y,z)
\end{equation*}%
and $r,r^{\prime },\beta $ are smooth functions of $z$ and $\delta $ is
constant. Further, $g(z)=\int \frac{1}{k_{3}(z)}dz.$

\begin{proof}
\textbf{\ }By virtue of (\ref{B7}) and (\ref{B8}), it can be proved that the
assumption $\nabla _{\xi }\tau =2a\tau \phi $ is equivalent to $\xi \cdot
\lambda =0.$ From the definition of gradient of a differentiable function we
get%
\begin{eqnarray}
grad\lambda &=&(e\cdot \lambda )e+(\phi e\cdot \lambda )\phi e+(\xi \cdot
\lambda )\xi  \label{grad} \\
&=&(e\cdot \lambda )e+(\phi e\cdot \lambda )\phi e  \notag
\end{eqnarray}%
Using (\ref{grad}) and $\left\Vert grad\text{ }\lambda \right\Vert =1$ we
have 
\begin{equation}
(e\cdot \lambda )^{2}+(\phi e\cdot \lambda )^{2}=1.  \label{1}
\end{equation}%
Differentiating (\ref{1}) with respect to $\xi $ and using$,$ (\ref{21}) and
(\ref{22}) we obtain,%
\begin{eqnarray*}
(\xi .e(\lambda )(e(\lambda )+(\xi \phi e(\lambda ))(\phi e(\lambda ) &=&0 \\
\left( \left[ \xi ,e\right] \left( \lambda \right) \right) e(\lambda
)+\left( \left[ \xi ,\phi e\right] \left( \lambda \right) \right) (\phi
e)\lambda &=&0\text{ \ \ } \\
\lambda e(\lambda )\phi e(\lambda ) &=&0
\end{eqnarray*}%
and since $\lambda \neq 0,$%
\begin{equation}
e(\lambda )\phi e(\lambda )=0  \label{25}
\end{equation}%
To study this system we consider the open subsets of U:%
\begin{equation*}
U^{\prime }=\left\{ p\in U\mid e(\lambda )(p)\neq 0\text{ }\right\} \text{,\
\ }U^{\prime \prime }=\left\{ p\in U\mid (\phi e)(\lambda )p\neq 0\right\}
\end{equation*}%
where $U^{\prime }\cup U^{\prime \prime }$ is open and dense in the closure
of $U.$ We distinguish two cases:

\textbf{Case1:}Now we suppose that $p\in U^{\prime }.$ By virtue of (\ref{1}%
), (\ref{25}) we have $(\phi e)(\lambda )=0,e(\lambda )=\mp 1.$ Changing to
the basis $(\xi ,-e,-\phi e)$ if necessary, we can assume that $e(\lambda
)=1.$ By the equation (\ref{22}), we get 
\begin{eqnarray}
\left[ \phi e,\xi \right] (\lambda ) &=&(\phi e)(\xi (\lambda ))-\xi ((\phi
e)(\lambda ))  \label{25a} \\
&=&(a-\lambda +1)e(\lambda )  \notag
\end{eqnarray}%
If we are using the relations $~e(\lambda )=1,(\phi e)(\lambda )=0$ and $\xi
\cdot \lambda =0$ in the equation (\ref{25a}), one can easily obtain $%
~a=\lambda -1.$ Hence, the equations (\ref{20}), (\ref{21}), (\ref{22}) and (%
\ref{3.5}), (\ref{3.6}) are reduced%
\begin{eqnarray}
\left[ e,\phi e\right] &=&-be+c\phi e+2\xi  \label{26} \\
\left[ e,\xi \right] &=&-2\lambda \phi e  \label{27} \\
\left[ \phi e,\xi \right] &=&0  \label{28}
\end{eqnarray}%
\begin{equation}
b=\frac{A}{2\lambda },\text{ \ }c=\frac{(B+1)}{2\lambda },~~  \label{29}
\end{equation}%
respectively.

Since $\left[ \phi e,\xi \right] =0,$ the distribution which is spanned by $%
\phi e$ and $\xi $ is integrable and so for any $p\in U^{\prime }$ there
exist a chart \{$V$,($x,y,z)\}$ at $p$, such that%
\begin{equation}
\xi =\frac{\partial }{\partial x}\text{, \ \ }\phi e=\frac{\partial }{%
\partial y},\text{ \ \ }e=k_{1}\frac{\partial }{\partial x}+k_{2}\frac{%
\partial }{\partial y}+k_{3}\frac{\partial }{\partial z}  \label{30}
\end{equation}%
where $k_{1},k_{2},k_{3}$ are smooth functions on $V$. Since $\xi $, $e$, $%
\phi e$ are linearly independent we have $k_{3}\neq 0$ $\ $at any point of $%
V $. Using (\ref{26}), (\ref{27}) and (\ref{30}) we get following partial
differential equations:%
\begin{equation}
\frac{\partial k_{1}}{\partial y}=\frac{A}{2\lambda }k_{1}-2\text{, \ \ }%
\frac{\partial k_{2}}{\partial y}=\frac{1}{2\lambda }\left[ Ak_{2}-B-1\right]
\text{ ,\ \ \ }\frac{\partial k_{3}}{\partial y}=\frac{A}{2\lambda }k_{3},
\label{31}
\end{equation}

\begin{equation}
\frac{\partial k_{1}}{\partial x}=0,\ \ \frac{\partial k_{2}}{\partial x}%
=2\lambda ,\ \ \frac{\partial k_{3}}{\partial x}=0,  \label{32}
\end{equation}%
Moreover we know that 
\begin{equation}
\frac{\partial \lambda }{\partial x}=0,\text{ \ \ }\frac{\partial \lambda }{%
\partial y}=0.  \label{33}
\end{equation}%
Differentiating the equation $\frac{\partial k_{3}}{\partial x}=0$ with
respect to $\frac{\partial }{\partial y}$ $,$ and using \ $\frac{\partial
k_{3}}{\partial y}=\frac{A}{2\lambda }k_{3}$ we find 
\begin{equation*}
0=\frac{\partial ^{2}k_{3}}{\partial y\partial x}=\frac{\partial ^{2}k_{3}}{%
\partial x\partial y}=\frac{1}{2\lambda }\frac{\partial A}{\partial x}k_{3}+%
\frac{1}{2\lambda }A\frac{\partial k_{3}}{\partial x}=\frac{1}{2\lambda }%
\frac{\partial A}{\partial x}k_{3}
\end{equation*}%
So

\begin{equation}
\frac{\partial A}{\partial x}=0  \label{34}
\end{equation}%
Differentiating $\frac{\partial k_{2}}{\partial x}=2\lambda $ with respect
to $\frac{\partial }{\partial y}$ $,$ and using $\frac{\partial k_{2}}{%
\partial y}=\frac{1}{2\lambda }\left[ Ak_{2}-B-1\right] $ and the equation (%
\ref{34}), we prove that

\begin{equation*}
\frac{\partial ^{2}k_{2}}{\partial y\partial x}=0=\frac{\partial ^{2}k_{2}}{%
\partial x\partial y}=\frac{1}{2\lambda }\left[ \frac{\partial A}{\partial x}%
k_{2}+A\frac{\partial k_{2}}{\partial x}-\frac{\partial B}{\partial x}\right]
\end{equation*}%
So%
\begin{equation}
\frac{\partial B}{\partial x}=2\lambda A  \label{35}
\end{equation}%
From (\ref{33}) we have following solution%
\begin{equation}
\lambda =\hat{g}(z)+d=\check{g}(z)  \label{36}
\end{equation}%
where $d$ is constant. Using $e(\lambda )=k_{1}\frac{\partial \lambda }{%
\partial x}+k_{2}\frac{\partial \lambda }{\partial y}+k_{3}\frac{\partial
\lambda }{\partial z}=1$ and (\ref{33}) we get 
\begin{equation}
\frac{\partial \lambda }{\partial z}=\frac{1}{k_{3}},\text{ \ }k_{3}\neq 0
\label{37}
\end{equation}%
If we differentiate the equation (\ref{37}) with respect to $\frac{\partial 
}{\partial y}$ and because of the equation $\frac{\partial \lambda }{%
\partial y}=0,$ we obtain%
\begin{equation}
0=\frac{\partial ^{2}\lambda }{\partial z\partial y}=\frac{\partial
^{2}\lambda }{\partial y\partial z}=-\frac{1}{k_{3}^{2}}\frac{\partial k_{3}%
}{\partial y}  \label{38}
\end{equation}%
Since $k_{3}\neq 0,$ the equation (\ref{38}) is reduced to%
\begin{equation}
\frac{\partial k_{3}}{\partial y}=0.  \label{39}
\end{equation}%
Combining (\ref{31}) and (\ref{39}) we deduced that 
\begin{equation}
A=0  \label{40}
\end{equation}%
Using (\ref{35}) and (\ref{40})we have%
\begin{equation}
\text{ }\frac{\partial B}{\partial x}=0  \label{41}
\end{equation}%
It follows from (\ref{41}) that 
\begin{equation}
B=F(y,z)  \label{42}
\end{equation}%
By virtue of (\ref{40}) , (\ref{31}) and (\ref{32}) we easily see that

\begin{equation}
k_{1}=-2y+r(z)  \label{43}
\end{equation}%
where $r(z)$ is integration function.

Combining (\ref{32}) and (\ref{39}) we get%
\begin{equation}
k_{3}=t(z)+\delta  \label{44}
\end{equation}%
where $\delta $ is constant.

If we use (\ref{32}), (\ref{36}), (\ref{40}) and (\ref{42}) in (\ref{31})%
\begin{equation}
\frac{\partial k_{2}}{\partial x}=2\check{g}(z),\ \ \frac{\partial k_{2}}{%
\partial y}=\frac{-(B+1)}{2\lambda }=\frac{-(F(y,z)+1)}{2\check{g}(z)}
\label{45}
\end{equation}%
It follows from this last partial differential equation that 
\begin{equation}
k_{2}=2x\check{g}(z)-\frac{(H(y,z)+y)}{2\check{g}(z)}+\beta (z)  \label{46}
\end{equation}%
where%
\begin{equation}
\frac{\partial H(y,z)}{\partial y}=F(y,z)  \label{47}
\end{equation}%
Because of (\ref{37}), there is a relation between $\lambda =\check{g}(z)$
and $k_{3}(z)$ such that $\check{g}(z)=\int \frac{1}{k_{3}(z)}dz$. We will
calculate the tensor fields $\eta $, $\phi $, $g$ with respect to the basis $%
\frac{\partial }{\partial x},\frac{\partial }{\partial y},\frac{\partial }{%
\partial z}$. For the components $g_{ij}$ of the Riemannian metric $g$,
using (\ref{30}) we have%
\begin{eqnarray*}
g_{11} &=&1,\text{ \ \ }g_{22}=1,\text{ \ \ }g_{12}=g_{21}=0,\text{ \ \ }%
g_{13}=g_{31}=\frac{-k_{1}}{k_{3}}, \\
g_{23} &=&g_{32}=\frac{-k_{2}}{k_{3}},\text{ \ \ \ \ \ \ }g_{33}=\frac{%
1+k_{1}^{2}+k_{2}^{2}}{k_{3}^{2}}.
\end{eqnarray*}%
The components of the tensor field $\phi $ are immediate consequences of%
\begin{eqnarray*}
\phi (\xi ) &=&\phi (\frac{\partial }{\partial x})=0,\text{ \ \ }\phi (\text{%
\ }\frac{\partial }{\partial y})=-k_{1}\frac{\partial }{\partial x}-k_{2}%
\frac{\partial }{\partial y}-k_{3}\frac{\partial }{\partial z}\text{\ ,} \\
\phi (\frac{\partial }{\partial z}) &=&\frac{k_{1}k_{2}}{k_{3}}\frac{%
\partial }{\partial x}+\frac{1+k_{2}^{2}}{k_{3}}\frac{\partial }{\partial y}%
+k_{2}\frac{\partial }{\partial z}.
\end{eqnarray*}%
The expression of the contact form $\eta ,$ immediately follows from%
\begin{equation*}
\eta =dx-\frac{k_{1}}{k_{3}}dz.
\end{equation*}%
Now we calculate the components of tensor field $h$ with respect to the
basis $\frac{\partial }{\partial x},\frac{\partial }{\partial y},\frac{%
\partial }{\partial z}.$%
\begin{eqnarray*}
h(\xi ) &=&h(\frac{\partial }{\partial x})=0,\text{ \ \ }h(\text{\ }\frac{%
\partial }{\partial y})=-\lambda \frac{\partial }{\partial y}\text{,} \\
h(\frac{\partial }{\partial z}) &=&\lambda \frac{k_{1}}{k_{3}}\frac{\partial 
}{\partial x}+2\lambda \frac{k_{2}}{k_{3}}\frac{\partial }{\partial y}%
+\lambda \frac{\partial }{\partial z}.
\end{eqnarray*}

\textbf{Case2:}Now we suppose that $p\in U^{^{\prime \prime }}.$ As in case
1, we can assume that $(\phi e)(\lambda )=1.$ Using the equations (\ref{20}%
), (\ref{21}) ,(\ref{22}) and (\ref{3.5}), (\ref{3.6}) are reduced%
\begin{eqnarray}
\left[ e,\phi e\right] &=&-be+c\phi e+2\xi  \label{48} \\
\left[ e,\xi \right] &=&0  \label{49} \\
\left[ \phi e,\xi \right] &=&-2\lambda e  \label{50}
\end{eqnarray}%
\begin{equation}
b=\frac{(A+1)}{2\lambda },\text{ \ \ }c=\frac{B}{2\lambda },~a=-1-\lambda
\label{51}
\end{equation}%
respectively.

Because of (\ref{49}) we find that there exist a chart \{$V^{\prime }$,($%
x,y,z)\}$ at $p$ $\in U^{^{\prime \prime }}$,%
\begin{equation}
\xi =\frac{\partial }{\partial x}\text{, \ \ }\varphi e=k_{1}^{\prime }\frac{%
\partial }{\partial x}+k_{2}^{\prime }\frac{\partial }{\partial y}%
+k_{3}^{\prime }\frac{\partial }{\partial z},\text{ \ \ }e=\frac{\partial }{%
\partial y}  \label{52}
\end{equation}%
where $k_{1}^{\prime },k_{2}^{\prime }$ and $k_{3}^{\prime }$ $%
(k_{3}^{\prime }\neq 0),$ are smooth functions on $V^{\prime }.$

Using (\ref{48}) ,(\ref{50}) and (\ref{52}) we get following partial
differential equations:%
\begin{equation}
\frac{\partial k_{1}^{\prime }}{\partial y}=\frac{B}{2\lambda }k_{1}^{\prime
}+2\text{, \ \ }\frac{\partial k_{2}^{\prime }}{\partial y}=\frac{1}{%
2\lambda }\left[ Bk_{2}^{\prime }-A-1\right] \text{ ,\ \ \ }\frac{\partial
k_{3}^{\prime }}{\partial y}=\frac{B}{2\lambda }k_{3}^{\prime },  \label{53}
\end{equation}

\begin{equation}
\frac{\partial k_{1}^{\prime }}{\partial x}=0,\ \ \frac{\partial
k_{2}^{\prime }}{\partial x}=2\lambda ,\ \ \frac{\partial k_{3}^{\prime }}{%
\partial x}=0,  \label{54}
\end{equation}%
Moreover we know that 
\begin{equation}
\frac{\partial \lambda }{\partial x}=0,\text{ \ \ }\frac{\partial \lambda }{%
\partial y}=0.  \label{55}
\end{equation}

As in Case 1, if we solve partial differential equations (\ref{53}) , (\ref%
{54}) and (\ref{55}) then we find \ 
\begin{equation}
B=0,\text{ \ \ }A=F(y,z)  \label{56}
\end{equation}%
\begin{equation}
\lambda =\bar{g}(z)+d^{\prime }=g(z),\text{ \ }k_{1}^{\prime }=2y+r^{\prime
}(z),\text{ \ }k_{3}^{\prime }=t^{\prime }(z)+\delta ^{\prime }  \label{57}
\end{equation}%
\begin{equation}
k_{2}^{\prime }=2xg(z)-\frac{(H(y,z)+y)}{2g(z)}+\beta ^{\prime }(z)
\label{58}
\end{equation}%
\begin{equation}
\frac{\partial H(y,z)}{\partial y}=F(y,z)  \label{59}
\end{equation}%
where $r^{\prime }(z)$ is integration function, $d^{\prime }$ and $\delta
^{\prime }$ are constants. By the help of (\ref{57}), the equation $(\phi
e)(\lambda )=1$ implies that 
\begin{equation}
\lambda (z)=g(z)=\int \frac{1}{k_{3}^{\prime }(z)}dz  \label{60}
\end{equation}

As Case1, we can directly calculate the tensor fileds $g,\phi $,$\eta $ and $%
h$ with respect to the basis $\frac{\partial }{\partial x},\frac{\partial }{%
\partial y},\frac{\partial }{\partial z}.$%
\begin{eqnarray*}
g &=&\left( 
\begin{array}{ccc}
1 & 0 & -\frac{k_{1}^{\prime }}{k_{3}^{\prime }} \\ 
0 & 1 & -\frac{k_{2}^{\prime }}{k_{3}^{\prime }} \\ 
-\frac{k_{1}^{\prime }}{k_{3}^{\prime }} & -\frac{k_{2}^{\prime }}{%
k_{3}^{\prime }} & \frac{1+k_{1}^{\prime 2}+k_{2}^{\prime 2}}{k_{3}^{\prime
2}}%
\end{array}%
\right) ,\text{ \ \ }\phi \text{\ }=\left( 
\begin{array}{ccc}
0 & k_{1}^{\prime } & -\frac{k_{1}^{\prime }k_{2}^{\prime }}{k_{3}^{\prime }}
\\ 
0 & k_{2}^{\prime } & -\frac{1+k_{2}^{\prime 2}}{k_{3}^{\prime }} \\ 
0 & k_{3}^{\prime } & -k_{2}^{\prime }%
\end{array}%
\right) , \\
\eta &=&dx-\frac{k_{1}^{\prime }}{k_{3}^{\prime }}dz\text{ \ \ and \ \ \ }h%
\text{\ }=\left( 
\begin{array}{ccc}
0 & 0 & -\lambda \frac{k_{1}^{\prime }}{k_{3}^{\prime }} \\ 
0 & \lambda & -2\lambda \frac{k_{2}^{\prime }}{k_{3}^{\prime }} \\ 
0 & 0 & -\lambda%
\end{array}%
\right) \text{\ }
\end{eqnarray*}
\end{proof}
\end{theorem}

\begin{remark}
By the help of Theorem 4, we can find a suitable coordinate chart of every
point of $M$ such that $A=0,B=F(y,z)$ or $A=F(y,z),B=0.$ If we put extra
assumption $F(y,z)=0$ relative to this chart then $\xi $ reeb vector field \
will be harmonic vector field. By the virtue the Theorem1.1 of \cite{KOU2},
the contact metric manifold $M^{3}$ becomes a generalized $(\kappa ,\mu )$%
-contact metric manifold.
\end{remark}

\begin{example}
We consider the 3-dimensional manifold%
\begin{equation*}
M=\{(x,y,z)\in R^{3},z\neq 0\}
\end{equation*}%
and the vector fields%
\begin{equation*}
\xi =\text{ }\frac{\partial }{\partial x},\text{ \ \ }\phi e=\frac{\partial 
}{\partial y},\text{ \ \ }e=-2y\frac{\partial }{\partial x}+(2xz-1)\frac{%
\partial }{\partial y}+\frac{\partial }{\partial z}.
\end{equation*}%
The 1-form $\eta =dx+2ydz$ defines a contact structure on $M$ with
characteristic vector field $\xi =\frac{\partial }{\partial x}$. Let $g$, $%
\varphi $ be the Riemannian metric and the $(1,1)$-tensor field given by 
\begin{eqnarray*}
g &=&\left( 
\begin{array}{ccc}
1 & 0 & a_{1} \\ 
0 & 1 & a_{2} \\ 
a_{1} & a_{2} & 1+a_{1}^{2}+a_{2}^{2}%
\end{array}%
\right) ,\text{ }\phi \text{\ }=\left( 
\begin{array}{ccc}
0 & a_{1} & a_{1}a_{2} \\ 
0 & a_{2} & a_{2}^{2}+1 \\ 
0 & -1 & -a_{2}%
\end{array}%
\right) , \\
\text{\ }h\text{\ } &=&\left( 
\begin{array}{ccc}
0 & 0 & -2yz \\ 
0 & -z & 2z(2xz-1) \\ 
0 & 0 & z%
\end{array}%
\right) ,\text{ \ \ }\lambda =z,
\end{eqnarray*}%
with respect to the basis $\frac{\partial }{\partial x},\frac{\partial }{%
\partial y},\frac{\partial }{\partial z}$, where $a_{1}=2y$ and $a_{2}=1-2xz$
.{\small \ }${\small \ }$By a straightforward calculation, we obtain 
\begin{equation*}
\nabla _{\xi }\tau =2(z-1)\tau \phi
\end{equation*}
\end{example}

Now, we will give an example satisfying Remark 2.

\begin{example}
We consider the 3-dimensional manifold%
\begin{equation*}
M=\{(x,y,z)\in R^{3},z>0\}
\end{equation*}%
and the vector fields 
\begin{equation*}
\xi =\text{ }\frac{\partial }{\partial x},\text{ \ \ }e=\frac{\partial }{%
\partial y},\text{ \ \ }\phi e=2y\frac{\partial }{\partial x}+(2xz-\frac{2z+y%
}{2z})\frac{\partial }{\partial y}+z\frac{\partial }{\partial z}.
\end{equation*}%
The 1-form $\eta =dx-\frac{2y}{z}dz$ defines a contact structure on $M$ with
characteristic vector field $\xi =\frac{\partial }{\partial x}$. Let $g$, $%
\varphi $ be the Riemannian metric and the $(1,1)$-tensor field given by 
\begin{eqnarray*}
g &=&\left( 
\begin{array}{ccc}
1 & 0 & -\frac{a_{1}}{a_{3}} \\ 
0 & 1 & -\frac{a_{2}}{a_{3}} \\ 
-\frac{a_{1}}{a_{3}} & -\frac{a_{2}}{a_{3}} & \frac{1+a_{1}^{2}+a_{2}^{2}}{%
a_{3}^{2}}%
\end{array}%
\right) ,\text{ \ \ }\phi \text{\ }=\left( 
\begin{array}{ccc}
0 & a_{1} & -\frac{a_{1}a_{2}}{a_{3}} \\ 
0 & a_{2} & -\frac{1+a_{2}^{2}}{a_{3}} \\ 
0 & a_{3} & -a_{2}%
\end{array}%
\right) , \\
\eta &=&dx-\frac{a_{1}}{a_{3}}dz\text{ \ \ and \ \ \ \ }h\text{\ }=\left( 
\begin{array}{ccc}
0 & 0 & -\lambda \frac{a_{1}}{a_{3}} \\ 
0 & \lambda & -2\lambda \frac{a_{2}}{a_{3}} \\ 
0 & 0 & -\lambda%
\end{array}%
\right) \text{\ }
\end{eqnarray*}%
with respect to the basis $\frac{\partial }{\partial x},\frac{\partial }{%
\partial y},\frac{\partial }{\partial z}$, where $a_{1}=2y,$ $a_{2}=2xz-%
\frac{2z+y}{2z}$, $a_{3}=z$ and $\lambda =\ln (z)$ .{\small \ }By direct
computations, we get 
\begin{equation*}
R(X,Y)\xi =(1-(\ln (z))^{2})(\eta (Y)X-\eta (X)Y)+2(-1-\ln (z))(\eta
(Y)hX-\eta (X)hY)
\end{equation*}%
and%
\begin{equation*}
\nabla _{\xi }\tau =2(-\ln (z)-1)\tau \phi
\end{equation*}
\end{example}

\end{document}